\input amstex
 \documentstyle{amsppt}
 \magnification=\magstep1
 \vsize=24.2true cm
 \hsize=15.3true cm
 \nopagenumbers\topskip=1truecm
 \headline={\tenrm\hfil\folio\hfil}

 \TagsOnRight

\hyphenation{auto-mor-phism auto-mor-phisms co-homo-log-i-cal co-homo-logy
co-homo-logous dual-izing pre-dual-izing geo-metric geo-metries geo-metry
half-space homeo-mor-phic homeo-mor-phism homeo-mor-phisms homo-log-i-cal
homo-logy homo-logous homo-mor-phism homo-mor-phisms hyper-plane hyper-planes
hyper-sur-face hyper-sur-faces idem-potent iso-mor-phism iso-mor-phisms
multi-plic-a-tion nil-potent poly-nomial priori rami-fication sin-gu-lar-ities
sub-vari-eties sub-vari-ety trans-form-a-tion trans-form-a-tions Castel-nuovo
Enri-ques Lo-ba-chev-sky Theo-rem Za-ni-chelli in-vo-lu-tion Na-ra-sim-han Bohr-Som-mer-feld}

\define\rest#1{_{\textstyle{|}#1}} 

\define\Span#1{\left<#1\right>} 

\define\id{\roman{id}} 

\define\C{\Bbb C} 
\define\R{\Bbb R} 
\define\Z{\Bbb Z} 

\define\proj{\Bbb P} 

\define\sA{{\Cal A}} 
\define\Oh{{\Cal O}} 

\define\al{\alpha}

\define\de{\delta}
\define\ep{\varepsilon}

\define\ga{\gamma}

\define\om{\omega}

\define\De{\Delta}




\redefine\mod{\mathop{{\roman{mod}}}} 


\define\Pic{\operatorname{Pic}} 



 \document

 \topmatter
  \title Three conjectures on lagrangian tori in the projective plane
 \endtitle
  \author    Nikolay A. Tyurin\footnote{This work was partially supported by RFBR, grants
  05 - 01 - 01086,  05 - 01 - 00455 and  07 - 01 - 92211.}
                \endauthor

   \address MPI (Bonn), BLTP JINR (Dubna)
   \endaddress
  \email ntyurin\@theor.jinr.ru
      jtyurin\@mpim-bonn.mpg.de
   \endemail

\abstract  In this paper we extend the discussion on Homological Mirror
Symmetry for Fano toric varieties presented in [HV] to more general case of
monotone symplectic manifolds with real polarizations. We claim that the Hori
-- Vafa prediction, proven in [CO] for toric Fano varieties, can be checked in
much more wider context. Then the notion of Bohr - Sommerfeld with respect to
the canonical class lagrangian submanifold appears and plays an important role.
The discussion presents a bridge between Geometric Quantization and Homological
Mirror Symmetry programmes both applied to the projective plane in terms of its
lagrangian geometry. Due to this relation one could exploit some standard facts
known in GQ to produce results in HMS.
\endabstract

   \endtopmatter

\head Introduction  \endhead

Lagrangian geometry of compact symplectic manifolds  remains to be a subject where not too much
is known. Even in the simplest case of two dimensional compact
symplectic manifolds (= Riemann surfaces) where the lagrangian condition
 degenerates and any 1- dimensional submanifold is lagrangian, the classification
problems (up to hamiltonian isotopy or up to symplectomorphism) are solved just
 for certain special cases including the case of the projective line. In dimension
4 one doesn't know which 2 - dimensional manifolds appear as lagrangian submanifolds,
 and the discussion on the existence of  a lagrangian Klein bottle
in  Kahler surfaces was not finished yet. One believes that the projective
plane $\C \proj^2$ admits only lagrangian tori as orientable lagrangian submanifolds
and real projective planes as non orientable ones plus some artificial types produced
 by hands using lagrangian surgery near the intersections of
lagrangian tori and real projective planes which gives new topological types of
 lagrangian submanifolds such as $T^2 \sharp \R \proj^2$ etc. On the other hand,
the classification   of lagrangian tori  even for $\C \proj^2$ is not completed:
one knows two types of lagrangian tori (the Clifford type and the Chekanov
type) and it seems that these types belong to different classes of the classification
 up to hamiltonian isotopy (one refers here to a paper of Chekanov and
Schlenk which is coming but  not published yet). And of course it doesn't mean
 that the set of equivalence classes is exhausted by these Clifford and Chekanov
types.

At the same time lagrangian geometry is highly desired for certain approaches to
Mirror Symmetry conjecture. According to  Homological Mirror Symmetry programme,
proposed by M. Kontsevich, [1], the symmetry is an equivalence between the
bounded derived category of coherent sheaves over a given algebraic manifold $M$
and certain Fukaya - Floer category of lagrangian submanifolds of its mirror
 partner $mir (M)$ where the last one is a symplectic manifold.
The objects of the Fukaya - Floer category are presented by lagrangian submanifolds
 up to hamiltonian isotopy, and the morphisms are given by
the Floer cohomology of pairs of objects. Despite of the fact that a rigouris
definition of the Floer cohomology doesn't exist in general,
one understands that the full category of lagrangian submanifolds up to hamiltonian
isotopy is too big anyway and this implies certain restrictions
on the type of the lagrangian submanifolds taken in the specified constructions.
 Concerning the 4 - dimensional case one adopted a variant of the generic Fukaya --
Floer theory where an additional data is exploited namely the structure of the
Lefschetz pencil on a given symplectic manifold. Then the Fukaya -- Seidel category
plays the role of the counterpartner of the bounded derived category of coherent
 sheaves; in the Fukaya -- Seidel category one takes not all lagrangian
submanifold but vanishing cycles only.

Another type of restriction is proposed for the case of toric Fano varieties: for
 such an $X$ one takes a toric fibration and considers the fibers $\{ S_{\al} \}$
(which are lagrangian submanifolds in $X$) with non trivial Floer cohomology
$FH^*(S_{\al}, S_{\al}; \C) \neq 0$. Then the desired category is constructed
over the set of the fibers, satisfy this non triviality condition. For this case
 the prediction of K. Hori and C. Vafa says that the number of such fibers is
finite; it should be the euler characteristic of the mirror partner, see [2]. This preduction
was proven in [3] modulo certain assumptions:instead of
the Floer cohomology one computes the Bott -- Morse version; instead of a generic
 almost complex structure one takes the complex structure of the toric
variety or its small hamiltonian deformations; the answer looks very familiar in
the framework of Geometric Qunatization of toric varieties -- the desired fibers
are distinguished by certain integrality condition, see formula (10.6) from [3].
 And this integrality condition in the toric framework means that one
deals with the Bohr -- Sommerfeld fibers.

But the main idea of Mirror Symmetry is to relate the algebraic geometry of a given
variety to the {\it symplectic geometry} of its mirror partner so
the answer on the right hand side must be independent on the choice of the compatible
 complex structure. One exploits certain sufficiently generic almost
complex structure to construct the objects like the Gromov invariants or the
Floer cohomology, and we know that an integrable complex structure is too far to be
generic in this setup (many examples of the answers which must be corrected are known
 in the gauge theories etc.); moreover the complex structure of a toric
variety is even more special (any toric variety is rigid in the class of toric varieties
 but of course can be deformed to a non toric algebraic variety).
Thus one needs to extend the setup of [3] in the way which would be more independent
on the complex structure choice. On the other hand, passing in this way
one can see that the results of [3] can be understood as more general fatcs adopted to
the specific situation of the toric Fano varieties.

What would be a possible "more general setup"? Instead of Fano varieties we consider
{\it monotone} simply connected compact symplectic manifold. Instead of
toric varieties we consider symplectic manifolds with {\it real polarizations}. And
after this translations we reach the situation which is
 very well known in Geometric Quantization, see f.e. [4], [5], [6].

 Geometric Quantization is a set of recepies attaching to a given symplectic manifold certain
Hilbert spaces together with homomorphisms of the Lie algebra of smooth functions on this
given manifold to the spaces of self adjoint operators acting on the
Hilbert spaces (or more generally, it attaches to a given symplectic manifold certain
 algebraic variety, [7]). One recepie from the set can be applied in the case when
our given symplectic manifold admits an additional structure --- a real polarization,
which is a lagrangian fibration of our given symplectic manifold.
In this case the Hilbert spaces are spanned by the fibers which satisfy some specific
 condition ---  so - called Bohr - Sommerfeld condition of different levels.
And the crucial fact here is the following: in the compact case the number of Bohr - Sommerfeld
 fibers is finite if the lagrangian fiberation is sufficiently good.
The toric Fano case is included by this class of "sufficiently good" lagrangian fibration
 and the construction from [3] ensures that for this case
there is the coincidence between Bohr -- Sommerfeld fibers and the fibers with non trivial
adopted Floer -- Bott -- Morse cohomology. To avoid the ambiguity
with different versions of the Floer cohomology we will universalize the story having in
 mind the following general distribution: if a lagrangian submanifold
is {\it displacable} then it seems that for any version of the Floer cohomology theory
it should have trivial Floer cohomology; and if a lagrangian submanifold is
{\it monotone} then it seems that for any definition it should have non trivial cohomology.
Recall that $S$ is {\it displacable} if it can be moved by
 some hamiltonian isotopy $\phi_H$ such that the intersection $S \cap \psi_H(S) = \emptyset$.
 Of course, the distribution is not complete in general,
but at least for the basic example of montone simply connected symplectic manifold
 with real polarization --- the projective plane ---
it seems to be exhaustive. The present text contains  results  about the fibers
of a real polarization with regular degeneration: for a generic monotone
symplectic manifold the number of monotone fibers is always finite (Theorem 2);
for any real polarization with regular degeneration of the projective plane any
Bohr -- Sommerfeld with respect to the canonical class fiber is monotone and
any other fiber is displacable (Theorem 3). Let us empasize that it is true for
{\it any} real polarization with regular degeneration, not just for toric one.
We do not use the toric structure in the constructions below, and in parallel
present several (naive) conjectures which extend the statement of Theorem 3 to
the case of {\it any} lagrangian torus in $\C \proj^2$, not just for the
fibers. One could not expect today that these ones are true but the work in
this direction is continued.

The discussion below follows the idea that  the Bohr - Sommerfeld condition and
the non triviality condition for the Floer cohomology are somehow related, and if it
 is indeed the case  we would get a way how to proceed  in HMS using known facts and
constructions
in GQ.  For the author the main reason to study the question is that it would be a
 realization of the idealogy, proposed by Andrey Tyurin, which claims that Mirror
 Symmetry and Geometric Quantization are relatives,  [8].

  I would like to thank the Max Planck Institute (Bonn) for the hospitality, excellent
 working condition and the possibilty to communicate with
  high level mathematicians. My gratitude goes to  F. Hirzebruch, A. Gorodentsev,
P. Pushkar, A. Kokotov, P. Moree and many others. I have to thank
  D. Auroux and D. Orlov for constant help and remarks. The last but not least thanks
 go to the staff of the Max - Planck - Institue fur Matematik.

  \head  \S 1. Bohr - Sommerfeld conditions \endhead

  Consider $(M, \om)$ --- a compact simply connected symplectic manifold of real dimension $2n$
and suppose that the cohomology
   class $[\om] \in H^2(M, \R)$ is integer.
  It means that there exists a complex line bundle $L \to M$ with the first Chern
class $c_1(L) = [\om]$. Choosing a hermitian structure on $L$,
  the space of hermitian connection $\sA_h(L)$ is defined. There exists unique up
 to gauge trnasformations hermitian connection $a \in \sA_h(L)$
  such that its curvature form $F_a$ is proportional to $\om$:
  $$
  F_a = 2 \pi i \om.
  $$
  The pair $(L, a)$ is usually called {\it the prequantization data}. Consider any
integer number $k \in \Z$ and the corresponding power $L^k$.
  The space $\sA_h(L^k)$ contains unique up to gauge transformations
  hermitian connection $a_k$ such that its curvature form is proportional to $\om$;
  $$
  F_{a_k} = 2 k \pi i \om.
  $$
  The pair $(L^k, a_k)$ is called the prequantization data {\it of level} $k$.

  Let $S \subset M$ be a lagrangian submanifold. This means that $S$ has dimension $n$ and
 the restriction $\om|_S$ vanishes identically.
  Then restricting the pair $(L^k, a_k)$ to $S$ one gets a trivial line bundle with a
flat connection since the curvature form vanishes being
  proportional to the symplectic form. Therefore lagrangian submanifolds can be distinguished
using the data of flat connections.

  We say that a lagrangian submanifold $S \subset M$ is {\it Bohr - Sommerfeld of level $k$}
 if the restricted connection $a_k|_S$ admits covariantly constant
  sections.

  One could ask whether or not this definition depends on the choices of hermitian structure
 on $L$ and a connection $a$ from the equivalence class of hermitian
  connections described by the condition $F_a = 2 \pi i \om$. The  point is that the
definition is {\it absolutely universal}: it can be reformulated in our case
  as follows. Consider $H_1(S, \Z)$ and for each primitive element $b \in H_1(S, \Z)$
consider some representative $\ga_b \subset S$. Since $\pi_1(M)$ is trivial
  one can find a disc $D \subset M$ with the boundary $\partial D = \ga_b$. It is not hard
to see that $S$ is Bohr - Sommerfeld of level $k$ if and only if
  for any $b \in H_1(S, \Z)$ the symplectic area of $D$ multiplied by $k$ is integer:
  $$
   k \cdot \int_D \om \in \Z.
  $$
  Note, that we consider the case when $[\om] \in H^2(M, \Z) \subset H^2(M, \R)$ and the
last integrality condition doesn't depend on the choice
  of particular $D$. At the same time the last description of the Bohr - Sommerfeld condition
doesn't involve any bundles or connections and therefore
  the notion is universal. It's natural to call the numbers
  $$
  p_k(b) = k \cdot \int_D \om  \quad mod \quad \Z
  $$
  {\it the periods} of our given lagrangian submanifold.

  On the other hand the last description relates the local deformations of a given symplectic
manfifold to the variations of  periods
  of the deformed submanifolds. According to the Darboux - Weinstein theorem, see [9], there
exists certain small tubular neighborhood
  of a given $S$ such that this neighborhood is symplectomorphic to a neighborhood of the
 zero section in $T^*S$, endowed with the standard
  symplectic form. Then local lagrangian deformations of $S$ are presented by graphs of closed
 1 -forms on $S$. Taking in mind the period
  description one can see that for a deformation $\psi$ of a given $S$ with the class $[\psi]
 \in H^1(S, \R)$ the periods change as follows
  $$
  p_k(\psi_*(b)) = p_k(b) + k \cdot \psi(b).
  \tag 1
  $$
  Indeed, if we have a loop $\ga_b$ on the given lagrangian submanifold and then deform it
to $\psi(\ga_b)$ the symplectic area of the tube with boundary
  $\ga_b - \psi(\ga_b)$ is exactly $\psi(b)$ and it implies the formula above.

  As a corollary one gets the existence of Bohr - Sommerfled lagrangian submanifolds of
(perhaps, sufficiently big) level $k$ in the case when a single
  lagrangian submanifold exists.

   A variant of the basic definition appears in the case of monotone symplectic manifolds.
A symplectic manifold $(M, \om)$ is called {\it monotone} if
   its canonical class $K_{\om} \in H^2(M, \Z)$ is proportional to $[\om]$:
   $$
   K_{\om} = k \cdot [\om].
   $$
   Number $k$ is called the coefficient of monotonicity. For this case for any hermitian
structure on $K_{\om}$ there exists unique up to gauge transformations
   hermitian connection $a_{can}$ with the curvature form proportional to the symplectic
form. Then we just repeat the basic definition with respect
   to $(K_{\om}, a_{can})$. If the restriction $(K_{\om}, a_{can})|_S$ admits covariantly
constant sections we say that $S$ is {\it Bohr - Sommerfeld with respect
   to the canonical class}. This specification is reasonable --- we will discuss it in the Section 3.

   \head \S 2. Finitness  \endhead

   Suppose now that we have an additional structure on $M$ --- a real polarization. This means
that $M$ is fibered over a base $B$ and almost all the fibers are
   smooth lagrangian. Usually it happens for the phase spaces of completely integrable
systems so there exists a commutative sub algebra in
   the Poisson algebra $(C^{\infty}(M, \R), \{, \}_{\om})$ spanned by the set of smooth
 functions $\{ f_1, ..., f_n \}$ called {\it integrals}
   such that the differentials $(d f_1, ..., d f_n)$ form a basis in the cotangent space
 almost everywhere on $M$. The conditions dictate several
   topological restrictions; the most important for us is that a smooth fiber must be
isomorphic to torus.

   For a real polarization of $M$ given by some map
   $$
   \pi: M \to B,
   $$
   where $B$ is a convex polytop in $\R^n$, we have the so - called {\it Kodaira --
Spencer map}; for each regular point $b \in B$ the deformation
    along the base is reflected somehow
   by the deformation of the fiber $\pi^{-1}(b) \subset M$ and since the local
deformations of a lagrangian submanifold are described by closed 1 - forms
   it induces a map
   $$
   m_{KS}: T_b B \to H^1(\pi^{-1}(b), \R).
   $$
   Now we claim that for a smooth fiber $\pi^{-1}(b)$ the Kodaira -- Spencer map is
 an  {\it isomorphism}. To prove this fact let us first mention that
   the dimensions of the entire spaces are the same:
   $$
   \dim T_b B = n = \dim H^1(T^n, \R).
   $$
   If one suppose that a vector $v \in T_b B$ goes to zero under the map then it would
imply that the corresponding small deformation is {\it isodrastic}
   so it preserves the periods of $S_b = \pi^{-1}(b)$. But according to  formula (1)
 it could happen if and only if the corresponding closed 1 - form
   $\psi$ is exact. This means that there exists a smooth function $f$ on $S_b$ such
 that $\psi = d f$. But each smooth function must have
   at least two critical points on a compact manifold, maximal and minimal. This means
 that the graph of $\psi$ in this case must intersect
   our given $S_b$ at least in two points. But it is impossible since the fibers cann't
 intersect each other. Therefore the Kodaira - Spencer map
   doesn't have a kernel and due to the dimensional reason it is an isomorphism.

    In Geometric Quantization the approach with real polarization gives the following
 recipe to construct the Hilbert spaces, see [5]. For a real polarization
    one takes the fibers which are Bohr - Sommerfeld of level $k$ which are
 $S_1, ..., S_l, ....$ and forms the linear span
    $$
    \sum_i \C <S_i > = {\Cal H}_k.
    $$
    The point is that the set of such fibers is discrete anyway and finite if the real
polarization has sufficiently good degenerations. Indeed, the discretness
    follows just from the fact that the Kodaira - Spencer map is an isomorphism. Now
 what are these "sufficiently good degenerations"? They appear
    for example in the case of toric varieties. This means that the degenerations are
 regular so if $B$ is a convex polytop in $\R^n$
    then the fibers over the inner part are smooth; the picture over  a $n-1$ -  dimensional
 face is again a smooth symplectic manifold fibered over
    this face with smooth lagrangian fibers over the inner part of this face (which
are smooth lagrangian $n-1$ - dimensional tori), etc. In this situation
    we have that the number of smooth Bohr - Sommerfeld lagrangian fibers is finite.
Indeed,  the discrete set can have a limiting point only on the boundary.
    Suppose that the limiting point corresponds to a smooth $n-1$ - dimensional torus
placed over the inner part of a $n-1$ - dimensional face. The preimage
    of this $n-1$ -dimensional face is a symplectic submanifold  $ M_1 \subset M$.
The limiting process implies  that for our fixed $k$ (the level
    of the Bohr - Sommerfeld property) the normal bundle $N_{M_1/M}|_{S_{lim}}$
contains a serie of shrinking dics bundles each of them consists of discs
    of constant symplectic area such that this area multiplied by $k$ is integer.
This implies that starting with some sufficiently small disc bundle
    the symplectic area of the fiber discs must be trivial. But it is impossible since
 the normal bundle $N_{M_1/ M}$ is symplectic and each disc must have
    nontrivial symplectic volume. Thus the limiting point on the inner part of a
$n-1$ - dimensional face cann't exist. Now suppose that the limiting point
    is more degenerated. In this situation one can use some natural shift of the
chain of lagrangian submanifold resulting with a generic limiting point which
    already lies on certain $n-1$ - dimensional face.

    Therefore we get the following

    {\bf Theorem 1.} {\it Let $X$ be a simply connected symplectic manifold and
 $\pi: X \to B$ be a real polarization with
    regular degeneration. Then for any level $k \in \Z$ the set of Bohr --
Sommerfeld lagrangian
    fibers of level $k$ is finite.}

   \head \S 3. Monotonicity \endhead

   As we've already mentioned in Section 1 there is a variant of the Bohr -- Sommerfeld
 condition natural in the setup of monotone symplectic manifolds.
   In the Floer cohomology theory one of the most important case is when a given
 lagrangian submanifold is {\it monotone}. For the specialization
   of such a manifold let us remind first what the Maslov index is. Since we are discussing
 below lagrangian tori  we consider the case of orientable
   lagrangian submanifolds.

   Let $S \subset M$ is an orientable lagrangian submanifold of a simply connected
 symplectic manifold $(M, \om)$. Choose any almost complex
   structure compatible with $\om$ and realize the anticanonical bundle $K^{-1}_{\om}$
 as the determinant of the hermitian bundle $(TM, I, \om)$.
   For any loop $\ga \subset S$ choose a disc $D \subset M$ with the boundary
 $\partial D = \ga$, and consider a trivialization of the anticanonical
   bundle  $K^{-1}_{\om}|_D$ restricted to $D$. This trivialization is unique up to
 gauge transformations and since $D$ is simply connected
   the degree of these transformations computed on the boundary must be trivial.
Due to the realization this trivialization is presented by a non
   vanishing on $D$ polyvector field $\eta$ of the type $(n, 0)$. On the other hand
 the boundary of the disc carries  a non vanishing real polyvector field
   $\theta$ which is given by the determinant of $TS$ restricted to $\ga$. Thus
 the hermitian pairing of $\eta$ and $\theta$ gives a map
   $$
   \phi_D: \ga \to \C^*,
   $$
   since it is not hard to see that the lagrangian condition implies that
 $<\eta, \theta>_h$ never vanishes. The degree of this map
   $$
   \mu(\ga, D) = <\phi^*_D h; [\ga]>
   $$
   (where $h$ is the generator of $H^1(\C^*, \Z)$ and $[\ga] \in H_1(\ga, \Z)$
 here is the fundamental class) is an integer number
   which doesn't depend on the choice of the almost complex structure. Moreover,
 it doesn't depend on the particular choice of $D$ in the same class
   from $\pi_2(M, S)$ with the image at $[\ga] \in \pi_1(S)$ under the canonical
 homomorphism. For another disc $D'$ with the same boundary $\ga$
    the value $\mu(\ga, D')$ can be computed in the following way:
    $$
    \mu(\ga, D') = \mu(\ga, D) + <K^{-1}_{\om}; [S^2 = D \cup D']>,
    $$
    hence if $D'$ is homotopy equivalent to $D$ then the numbers must be the same.
At the same time the number doesn't depend on the particular choice
    of $\ga$ in a given class $[\ga] \in \pi_1(S)$. Totally it shows that we have a map
    $$
    \mu: \pi_2(M, S) \to \Z,
    $$
    which is called the Maslov index. For any simply connected symplectic manifold and
 any lagrangian submanifold the index exists and moreover it is invariant
    under any lagrangian deformations. It easily follows from its definition --- it
 must be invariant under any contineous deformations. In the case when
    the ambient symplectic manifold has small second cohomology ($\Pic M = \Z$) the
 index can be reduced to a numerical correspondence
    $$
    \mu: H_1(S, \Z) \to \Z (\mod \deg K^{-1}_{\om})
    $$
    which is often called {\it the Maslov number}.

   Now, a lagrangian submanifold $S \subset M$ is  monotone if there exists an
 integer number $k$ such that for any loop $\ga \subset S$ and any disc $D \subset
   M, \partial D = \ga$ one has
   $$
   \mu(\ga, D) = k \cdot \int_D \om,
   \tag 2
   $$
   where $\mu$ is the Maslov index of the loop $\ga$ with respect to the disc $D$.
 The existence of a monotone lagrangian submanifold imposes strong
   restrictions on the topology of $M$ itself --- it must be monotone itself. And if
 it is monotone then it is reasonable to exploit the notion of
   Bohr - Sommerfeld with respect to canonical bundle lagrangian submanifolds.
 It is not hard to see that a lagrangian submanifold is monotone {\it only if}
   it is Bohr - Sommerfeld with respect to the canonical class. Indeed, the identity
 (2) is possible only in the case when for each $\ga, D$ the symplectic area
   of $D$, multiplied by $k$, is an integer number.   But it is exactly our Bohr
-- Sommerfeld condition with respect to the canonical class.

   On the other hand, for a Bohr - Sommerfeld with respect to the canonical bundle
 lagrangian submanifold in a monotone simply connected symplectic manifold
   one can define a characteristic class which is called {\it the universal Maslov class},
 see [10]. Leaving aside its first definition,
   we define it here as follows: for a given Bohr - Sommerfeld with respect to
 the canonical class $S \subset M$ with $K^{-1}_{\om} = k \cdot [\om]$
   for any loop $\ga$ and any disc $D, \partial D = \ga$ consider the difference:
   $$
   m_S(\ga, D) = \mu(\ga, D) - k \cdot \int_D \om \in \Z.
   $$
    Then the value of $m_S$ {\it doesn't depend} on the choice of $D$. Moreover,
 this numerical correspondence is {\it linear} and consequently
    $m_S$ is a {\it cohomology class} from $H^1(S, \Z)$. Let's remind that this
class is correctly defined if and only if our given lagrangian
    submanifold is Bohr - Sommerfeld with respect to the canonical bundle. And
 since this property is stable with respect to hamiltonian deformations only
    the resulting cohomology class is invariant under hamiltonian deformaitons
only, see [10].

    From this description we get tautologically that $S$ is monotone
{\it if and only if} it is Bohr - Sommerfeld with respect to the canonical class
    and its universal Maslov class vanishes
    $$
    m_S = 0.
    $$

    Resuming the discussion of this Section, we have the following proposition:

    {\bf Theorem 2.} {\it Let $X$ be a simply connected monotone symplectic manifold,
 and $\pi: X \to B$ be a real polarization with regular degeneration.
     Then the number of monotone lagrangian fibers is finite.}

    The prove is straightforward.

    \head \S 4. Lagrangian tori in $\C \proj^2$ \endhead

The resting part of the paper discusses the case of the basic example of the
monotone simply connected symplectic manifold --- the projective plane $\C
\proj^2$.

Take the projective plane $\C \proj^2$ with the standard Fubini --- Study
Kahler form $\om$
 which we consider as a symplectic form.
Thus the cohomology class $[\om]$ is integer and presents a generator
 of $H^2(\C \proj^2, \Z)$. As a symplectic manifold, it is monotone,
$K = -3 [\om]$.  We are interested in lagrangian fibrations of $\C \proj^2$
 to verify the following naive conjecture which can be attached to
any simply connected monotone symplectic manifold endowed with a real
 polarization with regular degeneration: if a smooth lagrangian fiber is
displacable it is not Bohr -- Sommerfeld with respect to the canincal class,
and if this fiber is Bohr -- Sommerfeld with respect to the canonical class it
is monotone (one could call this conjecture {\it  Extremely Naive}, or ENC for
short). Below we show that this conjecture is true for the projective plane.
But let us start with the basic example of lagrangian tori in $\C \proj^2$.

    The first and simplest example of lagrangian fibration of $\C \proj^2$
comes from the toric geometry;
    it can be given by the following construction. Choose  homogenious
 coordinates
$[z_0: z_1: z_2]$ and  consider a subset of $\C \proj^2$ defined by the
 system of equations
$$
z_i = r_i e^{ i \cdot\phi_i}, \quad i = 0, 1, 2,
$$
where $r_i$  are fixed positive real numbers satisfy $r_0 + r_1 + r_2 = 1$
and $\phi_i$ are real parameters. In $\C^3$ it would give a 3 -dimensional torus,
but after the phase factorization it gives us a smooth 2 - torus in $\C \proj^2$.
 Varying $r_i$s we get a family of lagrangian tori and hence a lagrangian
fibration of $\C \proj^2$ over a triangle $\De \subset \R^2$. Indeed, one can
 attach to any smooth torus the pair $(r_0, r_1)$ (since the third $r_2$
is defined by $(r_0, r_1)$ uniquelly), and the possible values of $(r_0, r_1)$
 form the triangle $\De$. The degenerations of this lagrangian fibration
are regular: over segments $\{ r_0 = 0, 0<r_1<1 \}, \{r_1 =0, 0<r_0<1\},
 \{r_0 + r_1 = 1, 0<r_i<1\}$ one has 1- dimenisonal torical fibers and the vertex
of the triangle $\De$ correspond to the maximal degenerations, 0 - dimensional
 tori or just points. Denote the smooth fiber of $(r_0, r_1)
\in \De$ as $S_{r_0, r_1}$. These are called {\it the Clifford tori} and the
 fibration is called {\it the Clifford fibration} of the projective plane.

Since the symplectic form is integer, the question arises about the Bohr --
Sommerfeld fibers of this lagrangian fibration. The line bundle $L = \Oh(1)$
with a hermitian connection $a$ whose curvature form is proportional to the
symplectic form distinguish a set of Bohr - Sommerfeld fibers of different
level. And the specification is very simple: the fiber $T_{r_0, r_1}$ is Bohr
-- Sommerfeld of level $k$ if and only if
$$
k \cdot r_0, k \cdot r_1 \in \Z.
$$
Indeed, the periods of the fiber torus $T_{r_0, r_1}$ are given by numbers
$r_0$ and $r_1$ for certain generators of $H_1(T_{r_0, r_1}, \Z)$
and this implies the statement.  This shows that:

--- there are no Bohr - Sommerfeld fibers of level 1 and 2;

--- there is unique fiber which is Bohr - Sommerfeld of level 3 and therefore
 which is Bohr - Sommerfeld with respect to {\it the canonical class};

--- the number of fibers which are Bohr - Sommerfeld of level $k$ is exactly
the same as $\dim H^0 (\C \proj^3, \Oh(k-3))$.

The last coincedence can be restored to direct equality "number of $k$ - Bohr
-- Sommerfeld fibers = dimension of holomorphic section space of $\Oh(k)$" if
one generalizes the situation and consider singular fibers as well. Then it
would be exactly three Bohr - Sommerfeld fibers of level 1 (= three points over
the vertices of $\De$, 0- dimensional tori), six Bohr - Sommerfeld fibers of
level 2 (= three points above + three middle 1- dimensional tori live over
edges of the triangle), ten Bohr - Sommerfeld fibers of level 3 (= three points
above + two for each edge of $\De$ 1 -dimensional tori + our regular fiber),
etc. This effect is known in Geometric Quantization of toric varieties.

But we are interested here in regular fibers only. Now let us see what would be
the result of the Homological Mirror Symmetry approach.
To proceed with one takes the fibers which have non trivial Floer cohomology.
 Leaving aside possible definitions of the Floer cohomology
$FH^*(S, S; \Z_2)$ we can exploit here our  remniscent: if a lagrangian
submanifold $S$ is {\it displacable} then it has trivial
Floer cohomology. And {\it the displacability} means that there exists a
 hamiltonian isotopy $\psi_t$ such that $\psi_t(S)$ doesn't intersect
$S$ for some $t$:
$$
\psi_t(S) \cap S = \emptyset.
$$

It is not hard to see that if both of $r_0$ and $r_1$ are not equal to $1/3$
then $T_{r_0, r_1}$ {\it is displacable}. Indeed, we have for $\C \proj^2$
the subalgebra of {\it symbols} in the Poisson algebra $(C^{\infty}(\C \proj^2),
 \{; \}_{\om})$, see [7], which correspond to self adjoint operators
on $\C^3$. The hamiltonian flow which moves $T_{r_0, r_1}$ to $T_{r_1, r_0}$
 is generated by the self adjoint operator which interchange
$z_0$ and $z_1$ in $\C^3$. And since the fibers don't intersect each other we
 get that if $r_0 \neq r_1$ then $T_{r_0, r_1}$ is displacable,
and the same is true if $r_1 \neq 1 - r_0 - r_1$. It remains one absolutely
 symmetric possibility: when all $r_i = 1/3$. And the point is that
this is precisely the Bohr -- Sommerfeld with respect to the canonical class
lagrangian fiber. To examine whether or not it has non trivial Floer cohomology
we use the following argument: the lagrangian torus $T_{1/3, 1/3}$ is monotone.
Indeed, it is Bohr - Sommerfeld with respect to the canonical class and it is
{\it minimal} therefore the universal Maslov class is trivial, see [10]. This
fact is exploited in [11] to prove that the Floer cohomology of $S$ is
isomorphic to the de Rham cohomology of it:
$$
FH^*(S, S; \C) = H^*_{dR}(S, \C),
$$
and the last one is very well known for a torus.

Thus for the standard toric fibration of $\C \proj^2$ (and the same is true
 for any projective space) the Bohr -- Sommerfeld with respect to the canonical class
condition is equivalent to the montonicity condition and furthermore to the non
 displacability condition (and in particular one could get
the results from [11]).

Now there is a natural simple extension of the toric case: a lagrangian torus
in $\C \proj^2$ is called
{\it of the Clifford type} if there exists a hamiltonian isotopy which moves
this torus to a standard fiber of the Clifford fibration. Since the Floer cohomology
 is invariant under hamiltonian deformations (and it is the main
property of it, which even could be taken for its general definition) as well
as the following three conditions:

--- the Bohr -- Sommerfeld condition of any level;

--- the monotonicity condition;

--- non displacability condition,

are, all what we've said is true for any Clifford torus.

Thus we complete the discussion of the Clifford tori in $\C \proj^2$, resuming that  for
the standard toric fibration of $\C \proj^2$ our
Extremely Naive Conjecture is true.

   The conjectures mentioned in the title of this text look rather naive as well being based
 mainly on the known examples and the facts that their statements are true  if we replace there
"a lagrangian torus" by "a fiber of a real polarization with regular
polarization", see below, but nevertheless we would like to formulte them in
these extended form.

   {\bf Conjecture 1.} {\it If $S \subset \C \proj^2$ is a Bohr -- Sommerfeld lagrangian
torus of level $k$ then $k$ must be greater of equal to
   3}.

   If this conjecture is true then the class of Bohr -- Sommerfeld with respect to
 the canonical class lagrangian tori is "pure" in the following sence:
   {\it a priori} any lagrangian torus which is Bohr - Sommerfeld of level 1 should
be automatically included to the set of Bohr - Sommerfeld
   with respect to the canonical class lagrangian tori (since 3 is divisible by 1),
but symplectically this torus is too far from
   the set of "pure" Bohr - Sommerfeld with respect to the canonical class tori.

   The next one is

   {\bf Conjecture 2.} {\it For any lagrangian torus $S \subset \C \proj^2$, Bohr -- Sommerfeld
 with respect to the canonical class,
   its universal Maslov class is trivial:}
   $$
   H^1(S, \Z) \ni m_S = 0.
   $$

   If this conjecture is true then every Bohr -- Sommerfeld with respect to
 the canonical class lagrangian torus should be monotone and thus must have
   non trivial Floer cohomology. This implies our third suggestion

   {\bf Conjecture 3.} {\it  A smooth lagrangian torus $S \subset \C \proj^2$ of
 the projective plane is non displacable} {\it if
    and only if $S$ is Bohr -- Sommerfeld with respect to the canonical class.}

   Let us note again that all these conjectures are too strong for proving a general
 version of the Hori -- Vafa prediction: it would be sufficient
    to exploit a weaker statement which looks as follows:

   {\bf Conjecture.} {\it Let $X$ be a monotone simply connected symplectic manifold
and $\pi: X \to B$ be a real polarization.  Then a smooth fiber
    $\pi^{-1}(b) = S_b$ is non displacable if it is Bohr - Sommerfeld with
respect to the canonical class.}

   Consider two examples both of non toric type.

   {\bf Toy example.} Consider $\C \proj^1 = S^2$ endowed with the standard symplectic
form. Any smooth loop $\ga \subset \C \proj^1$ is a lagrangian submanifold,
   and the topological type of smooth lagrangian submanifold actually is exhausted
by $T^1$, 1 - dimensional torus. Then the line bundle $L = \Oh(1)$ together with
   the appropiate hermitian connection $a \in \sA_h(L)$ defines the Bohr -- Sommerfeld
 condition of level $k$ which reads in this case as follows:
   a smooth loop $\ga \subset \C \proj^2$ is Bohr -- Sommerfeld of level $k$
if and only if it divides the surface into two pieces both of the symplectic area
   from $\frac{1}{k} \Z$. This means that $\ga \subset \C \proj^1$ is Bohr --
Sommerfeld with respect to the canonical class if and only if
   it divides the surface into  equal pieces. On the other hand, it is only
the case when $\ga$ is non displacable. This means that for a smooth loop in
   $\C \proj^1$ the Conjecture above is true.

  Of course, it says almost nothing for any other case since it is based on
the fact that for a smooth loop in $\C \proj^1$ there is only one symplectic
   invariant which characterizes the loop uniquelly up to symplectomorphism,
 namely the symplectic area of the disc, bounded by this loop.
   But for other dimensions it could be no longer true: one claims that there
 is at least one more type for lagrangian tori in $\C \proj^2$ which was called
  the lagrangian tori of the Chekanov type, see [12], [13]. Thus we have another

   {\bf "Non toric example".} This example can be found  in [13], where one
 characterizes it as a non toric fibration of $\C \proj^2$. It is defined as follows:
    consider the family of conics $\{ Q_{\ep} \}$ in $\C \proj^2$ given by the equation
    $$
    Q_{\ep} = \{ z_0 z_1 = \ep z_2 \},
    $$
    where $\ep \in {\bar \C}$, and $[z_0: z_1: z_2]$ is a homogenious coordinate
system. For this pencil with based points $[0: 0: 1], [0: 1: 0], [1: 0: 0]$
    one has exactly two singular conics:

    --- when $\ep = 0$ the conic is two intersecting lines;

    --- when $\ep = \infty$ the conic is   the double line $z_2^2 = 0$.

    Inside this pencil one could take 1- dimensional real subfamily consists of
the conics of the form $Q_{a \cdot e^{i t} - \mu}$ where $\mu \in \C^*,
    a \in \R^+$ are fixed numbers and $t$ is real parameter. For each entire of
this subfamily, say, $Q_{a \cdot e^{i t_0} - \mu}$ one has the natural fibration
    $$
    \pi: Q_{a \cdot e^{i t_0} - \mu} \to (-1, 1) \subset \R
    $$
    which is given by the hamiltonian action of the symbol which preserves each
 conics in the pencils (this simbol is essentially unique up to
    scale, see [14]). The fiber can be labeled by number $\de \in (-1, 1)$ and
 this correspondence has certain meaning: the symplectic area
    of the disc which bounds loop $T^1_{t_0, \de}$ equals to $\de \mod \Z$.
Fixing a value of $\de$ we distinguish the corresponding fiber $T^1_{t_0,\de}$.
    Now let us vary $t_0$ in the real subfamily $\{Q_{a \cdot e^{it} - \mu} \}$; it
gives us the corresponding family $T^1_{t, \de}$ which forms
    certain 2 -torus $T^2_{\de} = T^2_{a, \mu, \de} \subset \C \proj^2$. The point
 is that this torus is lagrangian, see [13]. All tori, constructed
    in this way for different $a$, form certain fibration of $\C \proj^2 \backslash l_2$
where $l_2$ is the line $z_2 = 0$ with only one singular torus
    $T^2_{\vert \mu \vert, \mu, 0}$ with one shrinked  loop (therefore strictly speaking
it is not a real polarization with regular degeneration).

    One claims that for a fixed $\mu \neq 0$ the fibration of $\C \proj^2$ consists
of two types of lagrangian tori:

    --- if $a > \vert \mu \vert$ then the torus $T_{a, \de}$ is of the Clifford type;

    --- if $a < \vert \mu \vert$ then the torus $T_{a, \de}$ is of the Chekanov type,

    see [13], and these types are different\footnote{It looks a bit strange since
 there is the case when $a = \vert \mu \vert$ and $\de \neq 0$.
    What is the type of this smooth lagrangian torus? It can be deformed to both
 the Clifford and the Chekanov types so if one takes in mind the
     fact that every contineous family of smooth lagrangian tori (submanifolds) with
 the same periods consists of hamiltonically equivalent tori,
     the types should be related.}.

    What it gives for our discussion and conjectures? Let us note that

    --- every torus $T^2_{a, \de}$ of the Clifford type has been discussed above
 so the conjectures are true for them;

    --- every torus of Chekanov type $T^2_{a, \de}, a < \vert \mu \vert,$ is
 {\it displacable} and thus has trivial Floer cohomology to itself;

    --- there is no a torus of Chekanov type which is Bohr -- Sommerfeld with
respect to the canonical class.

    Indeed, it is not hard to construct a smooth function on $\C \proj^2$ whose
hamiltonian flow moves $T_{a, \mu, \de}^2$ to $T^2_{a, - \mu, \de}$.
    This function is not generic, it is a symbol which corresponds to the
self adjoint operator $A = {\roman diag} (0, 0, 1)$. The hamiltonian flow
    then acts as the rotation of the parameter space $\C$ with two fixed points
$0$ and $\infty$. Then since a torus of the Chekanov type corresponds
    to the case when two circles of the same radius $a$ with centers at
$\mu$ and $- \mu$ do not intersect each other one gets that it is displacable.
    On the other hand,   direct computations show that there is no Bohr -- Sommerfeld
with respect to the canonical class lagrangian fiber of the Chekanov
    type.

    However the last example is excluded by our main setup of real polarizations with
regular degeneration. Regularity of degenerations
    took place for a given case of a real polarization imposes a number of natural arguments and facts.

    Consider any real polarization  $\pi: \C \proj^2 \to B \subset \R^2$ where $B$ is a
 convex polytop. Suppose that it has only regular degenerations. This means
    that there exists a set of symplectic divisors $D_1, ..., D_m \subset \C \proj^2$ such
 that $\dim_{\R} D_i = 2$. These divisors lie over
    the edges of $B$. Then it follows that

    --- each $D_i$ represents the class $[D] \in H_2(\C \proj^2, \Z)$ Poincare dual to
the cohomology class $[\om]$;

    --- the number of the symplectic divisor is $\deg K^{-1} = 3$.

    Indeed, the total degree of the boundary components must be the degree of the anticanonical
 class since the "inner" part of $\C \proj^2$ modulo
    $B$ admits non vanishing holomorphic vector 2 - field with respect to an almost
 complex structure, compatible with $\om$ and $\pi$. Each component
    from $\pi^{-1}(\partial B)$ must have positive degree with respect to $[\om]
 \in H^2(\C \proj^2, \Z)$. On the other hand the number of components
    equals to the number of edges of our convex polytop $B$ which is must be greater
 or equal to 3. This shows that for $\C \proj^2$ regularity dictates
    the form of $B$ and the type of $\pi^{-1}(\partial B)$.

    Furthemore, the analysis of the system would become simpler if it were possible
 to find integrals of special type. Namely, since the "inner"
    part  of $\C \proj^2$ is topologically equivalent to the direct product
$(B - \partial B) \times T^2$ we can choose a basis in
    $H_1(\pi^{-1}(b), \Z)$ uniformally for all smooth fibers of $\pi$. Moreover
it can be done relatively to the boundary components $D_1, D_2, D_3$
    if one chooses any two from the set. The point is that for any $D_i$ there
exists uniquelly determined basic element from $H_1(\pi^{-1}(b), \Z)$
    which degenerates when passing to a limit fiber in $D_i$. This means that we
 have distinguished primitive elements $d_1, d_2, d_3 \in H_1(\pi^{-1}(b), \Z)$,
    such that each pair $d_i, d_j, i \neq j,$ form a basis. Let us choose and fix $d_1, d_2$
 as a basis. Then there exists a lift of the period map with
    respect to the boundary data:
    $$
    p_{d_1, d_2} = (p^1_{d_1, d_2}, p^2_{d_1, d_2}): B \to \R^2,
    $$
    such that

    ---   $p_{d_1, d_2}$ is smooth on $B - \partial B$;

    --- $p^i_{d_1, d_2}|_{\pi(D_i)} = 0$;

    --- $p^i_{d_1, d_2} (b) = \int_D \om \mod \Z$ where $D \subset \C \proj^2$ is a disc
 with boundary $\partial D = \ga_{d_i} \subset \pi^{-1}(b)$
    and $[\ga_{d_i}] = d_i \in H_1(\pi^{-1}(b), \Z)$.

    Note that such a lift exists and there are exactly 4 possibilities for the extension $p_{d_1, d_2}$
 since there are exactly 4 possible choices of the signs for $d_1$ and $d_2$
    (compare this fact with the discussion on the choice of  spin structures in [3], [11]).
 Let us fix the signs in such a way that
    $p_{d_1, d_2}$ is non negative on $B$.

    Denote as $a_{ij}$ the intersection points, such that
    $$
    a_{ij} = D_i \cap D_j.
    $$
    Then it is easy to see that
    $$
    p_{d_1, d_2} (a_{12}) = (0, 0),  p_{d_1, d_2}(a_{13}) = (0,1),  p_{d_1, d_2}(a_{23}) = (1, 0).
    $$
    From this one deduce that
    $$
    d_3 = d_1 + d_2 \in H_1(\pi^{-1}(b), \Z).
    $$
    Indeed, $d_3$ can be represented as $ d_3 = p d_1 + q d_2$ being a primitive element,
where $p, q$ are coprime integers. But the symplectic area is an additive
    functional and from this one deduces that $p = q = 1$.

    Now let us impose the fact, proven in Section 2: the Kodaira -- Spencer map is
 an isomorphism. This implies one very important property of our lifted period function:

    {\bf Lemma.} {\it The function $p_{d_1, d_2}$ is strictly monotone in both arguments.}

Indeed, since each component of $p_{d_1, d_2}$ is monotone on the corresponding
boundary side and the fact, that the Kodaira - Spencer map in this situation
coincides with {\it the differential} of $p_{d_1, d_2}$ one sees that

--- the lifted period map $p_{d_1, d_2}$ doesn't have any critical points on $B - \partial B$;

--- for any level line $ L_c = \{ p^i_{d_1, d_2} = c, 0 \geq c < 1 \}$ the restriction
$p^j_{d_1, d_2}|_{L_c}$ is a strictly monotone (increasing) function.

Now examine the statements of Conjectures 1 -- 3 for this situation.

{\bf Conjecture 1.} From the monotonicity of $p_{d_1, d_2}$ it follows that for
any regular fiber $S_b = \pi^{-1}(b), b \in B - \partial B$, one has
$$
 0 < p^1_{d_1, d_2}(b) + p^2_{d_1, d_2}(b) < 1.
$$
By the definitions of $p_{d_1, d_2}$ and of the Bohr -- Sommerfeld fiber of
level $k$ we get that the minimal possible non empty level is 3.

{\bf Conjecture 2.} Again from the monotonicity of $p_{d_1, d_2}$ we get that
there exists unique fiber which is Bohr -- Sommerfeld of level 3 or with
respect to the canonical class. Note that for this fiber $S_{can}$ one has
$$
p^1_{d_1, d_2}(S_{can}) = p^2_{d_1, d_2}(S_{can}) = \frac{1}{3}.
$$
To prove the monotonicity of $S_{can}$ it is sufficient to find for each
generator of $H_1(S_{can}, \Z)$ a smooth loop $\ga$ representing this
generator, and a smooth disc $D$, bounded by $\ga$, such that the Maslov index
of $[\ga, D]$ would be three times the symplectic area of $D$ (and it is enough
since for any other disc $D'$ with the same boundary $\ga$ the relation should
be the same due to the monotonicity of $\C \proj^2$). Note that since the set
of lagragnian fibers is connected the Maslov index is the same for all
lagrangian tori.

 Take our distinguished  generator $d_1 \in H_1(S_{can}, \Z)$ and choose a
smooth loop $\ga_1 \subset S_{can}$ such that $[\ga_1] = d_1$. Take the level
set $C_{\frac{1}{3}} = \{ p^2_{d_1, d_2} = \frac{1}{3} \}$ and choose the
segment $ B_t \subset C_{\frac{1}{3}}, t \in [0; \frac{1}{3}],$ which
corresponds to the inequality $p^1_{d_1, d_2} \leq \frac{1}{3}$. There exists a
family of smooth loops $\ga_1^t, t \in [0; \frac{1}{3}]$ such that

--- $\ga_1^{\frac{1}{3}} = \ga_1 \subset S_{can}$;

--- $\pi(\ga_1^t) = b(t) \in B_t \subset B$;

--- $\ga_1^t \subset S_{b(t)}$ and $[\ga_1^t] = d_1 \in H_1(S_{b(t)})$.

The point is that the family $\{ \ga_1^t \}$ shrinks to point $\ga_1^0$ which
lies on the symplectic divisor $D_1$.

It is not hard to see that the family $\{ \ga_1^t \}$ forms a disc
$$
\cup_{t \in [0; 1/3]} \ga_1^t = D \subset \C \proj^2
$$
 such that
$$
\int_D \om = \frac{1}{3}.
$$
On the other hand the maslov index of $[\ga_1, D]$ is equal to 1. Indeed, since
we shrink $\ga_1$ to a point over the level line of $p^2_{d_1, d_2}$ it follows
that the Maslov index of $D$ must be the degree of the normal bundle of $D_1$.
Thus  we have
$$
\mu([\ga_1, D]) = 1 = 3 \cdot \frac{1}{3} = 3 \cdot \int_D \om,
$$
and since we can repeat the arguments for a smooth loop $\ga_2 \subset
S_{can}$, which represents the generator $d_2$, and it follows that $S_{can}$
is monotone.

{\bf Conjecture 3.} To prove the fact that if a fiber $S_b$ is not Bohr --
Sommerfeld with respect to the canonical class then it is displacable it is
sufficient to prove that the same happens for the fiber which doesn't lie over
the "diagonal" $\{ p^1_{d_1, d_2} = p^2_{d_1, d_2} \}$. Indeed, our choice of
$d_1, d_2$ was made arbitrary and taking another pair, say, $(d_1, d_3)$ we
shall get the same result for the corresponding "diagonal", and since the
intersection of the "diagonals" consists of exactly one point which is Bohr --
Sommerfeld with respect to the canonical class, it implies that the Conjecture
3 is true for fibers of a real polarization with regular degeneration of $\C
\proj^2$.

We calim that there exists a Hamiltonian deformation of $\C \proj^2$ which
generates the corresponding Hamiltonian isotopy which interchanges fibers with
values $(c_1, c_2)$ and $(c_2, c_1)$ with respect to the function $p^i_{d_1,
d_2}$. The desired Hamiltonian deformation is constructed explicitly as
follows. Consider the level sets of the sum $p^1_{d_1, d_2} + p^2_{d_1, d_2}$,
lifted to $\C \proj^2$ . The possible values are in $[0;1]$. There are two
exeptional level sets: for $c = 0$ we have the point $D_1 \cap D_2$; for $c =
1$ it is $D_3$. For any other $\al \in (0; 1)$ the level set
$$
C_{\al} = \pi^{-1}(\{ p^1_{d_1, d_2} + p^2_{d_1, d_2} = \al \}
$$
is a smooth 3 - sphere. The restriction of the symplectic form $\om$ to
$C_{\al}$ defines a fibration
$$
p_{\al}: C_{\al} \to S^2_{\al}
$$
which is topologically the Hopf bundle. Indeed, we take the kernels of
$\om|_{C_{\al}}$, and the corresponding 1 - dimensional distribution is
integrable which gives the fibration. Additionally one has

--- a symplectic form $\om_{\al}$ on $S^2_{\al}$ which is the result of the reduction applied to $\om$;

--- a smooth circle $S^1_{\al} \S^2_{\al}$ which is the result of the phase factorization
of the "diagonal" torus with periods $(\al/2, \al/2)$.

Note that when $\al$ tends to 1 this Hopf bundle $C_{\al} \to S^2_{\al}$
degenerates to $D_3$ with a marked circle $S^1_1 \subset D_3$. Moreover, the
triple $(D_3, \om|_{D_3}, S^1_1)$ is the result of the limiting procedure
applied to $(S^2_{\al}, \om_{\al}, S^1_{\al})$ when $\al$ tends to 1. On the
other hand, the other limit $\al \to 0$ is realized as a conformal shrinking of
the triple $(S^2_{\al}, \om_{\al}, S^1_{\al}$ to the point $D_1 \cap D_2$.
Indeed, it is clear that the symplectic volume
$$
\int_{S^2_{\al}} \om_{\al} = \al.
$$

 Let us fix for the symplectic 2 - sphere $D_3$ a smooth function $f_1 \in C^{\infty}(D_3, \R)$
such that $f_1$ is a height function and it has two non degenerated critical
points $p_1^N, p_1^S$ both of which lie on the marked circle $S^1_1$. One can
consruct now using inverse limiting process a family of smooth functions $\{
f_{\al} \}$ for the family of 2 - spheres $\{ S^2_{\al} \}$ for $\al \in (0;
1]$. We take an appropriate normalization for the functions such that
$$
\int_{S^2_{\al}} f_{\al} \om_{\al} = \al^2,
$$
and then lift each function $f_{\al}$ to $C_{\al}$ via the canonical
projection:
$$
F_{\al} = f_{\al} \circ p_{\al}: C_{\al} \to \R.
$$
  Then we claim that these lifted functions can be combined to a global smooth
function $F$ such that
$$
F|_{C_{\al}} = F_{\al}.
$$
This function has exactly three critical points:

--- the intersection point $D_1 \cap D_2$,

--- two points $p_1^N, p_1^S$ which lie on $D_3$.

The hamiltonian vector field $X_F$ generates the flow on $\C \proj^2$ which is
a 1 - parameter family of symplectomorphisms $\phi_t$ of $\C \proj^2$ such that
$\phi_{2\pi} = \id$. Indeed, it follows the rotation of the spheres $S^2_{\al}$
with fixed points $p^N_{\al}, p^S_{\al}$ which lie on the "diagonal" circle
$S^1_{\al}$. And it is not hard to see that the result of this rotation applied
to  a fiber of the given real polarization with periods $c_1, c_2)$ should be
the fiber with periods $(c_2, c_1)$. It ends the prove of the Conjecture 3 for
fibers of real polarization with regular degeneration.

Resuming the discussion we see that the following fact takes place:

{\bf Theorem 3.} {\it For $\C \proj^2$ the Extremely Naive Conjecture is true}.

Note that the method of lifting of the period map can be applied to any compact
simply connected symplectic manifold, and the main property of the Kodaira --
Spencer map can be exploited to establish the strict monotonicity of this
lifted period function which was crucial in our construction for $\C \proj^2$
above. Thus one could expect that the same method will be usefull for more
general cases, for other monotone symplectic manifolds.

At the same time before studing the Conjectures 1 -- 3, which were formulated
for any lagrangian tori in $\C \proj^2$, one could try to find the answer on
the following natural question: is there a geometric condition on a lagrangian
torus in $\C \proj^2$ which should detect whether or not this torus can be
included to a family of lagrangian fibers of a real polarization with regular
degeneration?

\enddocument